# On Hadwiger's Conjecture

Dhananjay P. Mehendale
Sir Parashurambhau College, Tilak Road, Pune 411030,
India

## Abstract

We propose an algorithm to reduce a *k*-chromatic graph to a complete graph of largest possible order through a well defined sequence of contractions. We introduce a new matrix called transparency matrix and state its properties. We then define correct contraction procedure to be executed to get largest possible complete graph from given connected graph. We finally give a characterization for *k*-chromatic graph and use it to settle Hadwigers conjecture.

**1. Introduction:** In the year 1943 Hadwiger proposed the following conjecture [1] which states that "Every *n*-chromatic graph contains a subgraph contractible to $K_n$, a complete graph on *n* points". For *n* = 5 this conjecture is equivalent to the famous "Four color problem".

The four color problem first appeared in a letter of October 23, 1852 to Sir William Hamilton from Augustus de Morgan, which was asked to him by his student Frederick Guthrie who later attributed it to his brother Francis Guthrie.

After the announcement of this problem to the London Mathematical Society by Arthur Cayley in 1878, within a year its solution was proposed by Kempe [2]. One year after this publication Heawood published its refutation [3]. But Kemp's idea of alternating paths actually led to the correct solution of this problem in the hands of Appel and Haken, obtained by making heavy use of computer [4]. Before we proceed to discuss the notion of transparency matrix and algorithm that leads to a complete graph of largest size through a sequence of well defined contraction operations which we think will be very useful for Hadwiger's problem a brief description of some useful ideas is in order. We will discuss them in the next section. We have followed [5] for their description.

**2. Preliminaries:** Let $G$ be a graph on $p$ points (vertices) and $q$ lines (edges) with vertex set $V(G) = \{v_1, v_2, \cdots, v_p\}$ and edge set $E(G) = \{e_1, e_2, \cdots, e_q\}$. A $k$-coloring of a graph $G$ is a labeling defined by the mapping $f : V(G) \to \{1, 2, \cdots, k\}$. The set vertices with same color form a color class. A $k$-coloring is proper if $x$ is adjacent to $y$ implies $f(x) \neq f(y)$. A graph is $k$-colorable if it has a proper $k$-coloring. The chromatic number, $\chi(G)$, is the minimum $k$ such that $G$ is $k$-colorable. If $\chi(G) = k$, then $G$ is $k$-chromatic. The edge contraction is the operation of putting vertex $x$ on vertex $y$, or vice versa, and is defined when vertex $x$ is adjacent to vertex $y$. $G$ is **contraction sensitive** if $\chi(G) = k$, but $\chi(H) < k$ for every graph $H$ obtained by the contraction of any edge of $G$. A connected graph is called $k$-partite if its vertex set can be represented as union of totally disconnected sets, $k$ in number, such that each of these sets form a color class (independent sets) and all the edges of the graph have end vertices in some different sets among these $k$ sets. Every connected graph $G$ can be looked upon as a $k$-partite graph for some $k$. The representation of $G$ as some $k$-partite graph where $k$ is minimum (i.e. $G$ does not have a representation as an $l$-partite graph where $l < k$) is called the **minimal-partite-representation** of $G$. For every connected graph there will exist a minimal-partite-representation. A partite set among the partite sets, $k$ in number, a set containing $m$ elements, $m > 1$, is called **essentially singleton** if there exists another partite set containing $n$ elements in which all the elements of earlier set can be added, except some one element, and still both the newly formed partite sets, along with other partite sets, remain independent. Every graph as some $k$-partite graph is clearly $k$-colorable and it will be $k$-chromatic when $k$ is minimum and every partite set forming a color class. The alternative definition of a contraction sensitive graph can be as follows: A connected $k$-chromatic graph is **contraction sensitive** if it has a minimal-partite-representation as a $k$-partite graph but if **any one** of its edge is contracted then it has a minimal-partite-representation as a $l$-partite graph such that $l < k$. A pair of vertices in a partite set is called **separators** if these vertices are simultaneously adjacent to some vertex in the other partite set. Thus, the separators are separated from each other by distance two (a two-path). It is clear that when some one edge in this two path is contracted then such contraction splits the partite set in which the separators belong into two partite sets.

**3. Hadwiger's Conjecture:** Let $G$ be a $k$-chromatic graph on $p$ points and $q$ lines, $p > k$. Clearly this graph will always have a minimal-partite-representation, as some $k$-partite graph. Because, if $k$ is not minimum and suppose there exists subpartite-representation for $G$ as some $l$-partite graph with $l < k$, then $G$ will be $l$-chromatic, a contradiction.

In order to settle Hadwiger's problem we essentially need to show that for a $k$-chromatic graph $G$ there exists a well defined sequence of contraction operations to be carried out which will take it to $K_k$, or to a graph which contains $K_k$ as its subgraph. In other words, we need to show that a minimal-partite-representation for a $k$-partite graph does not remain contraction sensitive at each stage of contraction and on the contrary after each contraction sensitive stage the separators will get crated which will restore of its original $k$-chromaticity. We characterize below contraction sensitive graph:

**Theorem 3.1:** A $k$-chromatic graph on $p$ points and $q$ lines, $p > k$, is contraction sensitive if and only if every partite set is essentially singleton.

**Proof:** Suppose the figure below represents the minimal-partite-representation of $G$ as $k$-partite graph. Thus, $V(G) = \bigcup_{i=1}^{k} A_i$, such that $A_i \cap A_j = \phi$, a null set when $i \neq j$, for all $i, j = 1, 2, \cdots, k$ and $A_i$ form the independent sets. Note that in this Figure 1 we have shown only the representative edges connecting the $k$ partite sets.

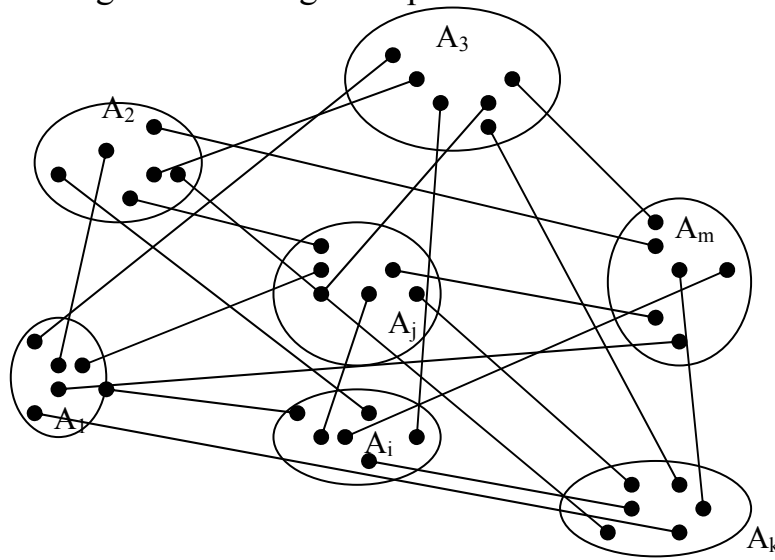

**Figure 1**

Suppose there exists no edge which when contracted produces a graph on (*p*-1) points and (*q*-1) lines which is still a *k*-partite graph with *k* minimum, i.e. this representation is contraction sensitive and if we contract any of its edge then it results in the reduction in the count of partite sets to (*k*-1), i.e. let there exists a vertex $x_i$ in the partite set $A_i$ adjacent to vertex $x_j$ in the partite set $A_j$ such that when the edge joining $x_i$ and $x_j$ is contracted the sets $A_i$ and $A_j$ get merged into each other to form a new independent partite set $A_i \cup A_j$ with $x_i$ and $x_j$ are now identified. But, in such case the original partitioning of the vertex set given above can be represented in the following **alternative** way:

$A_1 = A_1$, $A_2 = A_2$, $A_3 = A_3, \cdots, A_i = \{x_i\}$, $A_j = A_i \cup A_j - \{x_i\}$, $A_{i+1} = A_{i+1}, \cdots, A_k = A_k$.

This partitioning of the vertex set into partite sets will produce the following Figure 2.

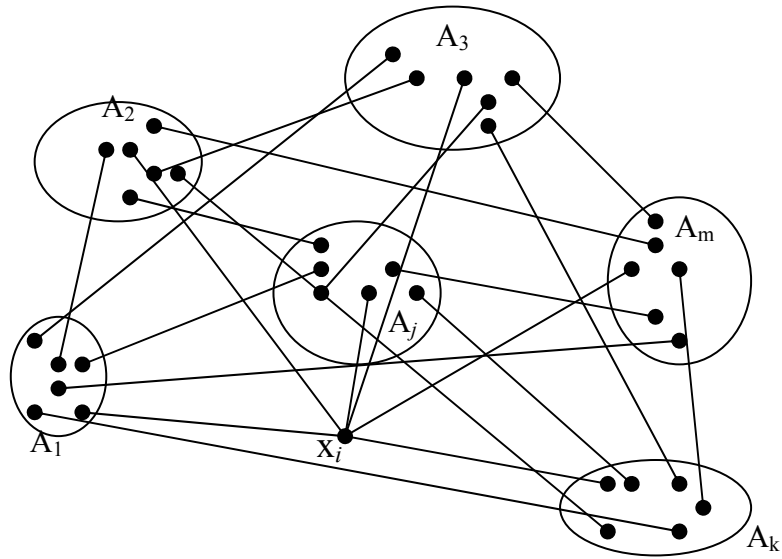

**Figure 2**

Thus, the partite set $A_i$ is essentially singleton.

Converse is straightforward. □

**4. Transparency Matrix:** Now, the immediate question that arises is the following: How one can search, in a simplest possible way, for a minimal-partite-representation for given graph? We now define a new matrix called **transparency matrix** and study its properties. We will see that these properties can be very useful for one's search for a proof of the Hadwiger's conjecture.

**Definition 4.1: Transparency matrix**, $T(G)$, associated with a graph $G$ containing $p$ points and $q$ lines, is the following $p \times p$ matrix:
$$T(G) = [a_{ij}]_{p \times p}$$
where, $a_{ii} = 0$, $a_{ij} = k$, where $k$ is the distance between vertices $v_i$ and $v_j$, i.e. it is the length of the shortest path joining vertices $v_i$ and $v_j$.

We now note certain **interesting properties** of $T(G)$:

(1) When there is no path connecting vertices $v_i$ and $v_j$, i.e. when $G$ is a disconnected graph and $v_i$ and $v_j$ belong to two different connected components of $G$, then $a_{ij} = \infty$. Also,

(2) When $v_i$ and $v_j$ are adjacent vertices then $a_{ij} = 1$.

(3) If we replace all the so called distances $k \geq 2$ by zero then the transparency matrix becomes the usual adjacency matrix, $A(G)$.

(4) All the two by two principle sub-matrices of the form $\begin{bmatrix} 0 & 1 \\ 1 & 0 \end{bmatrix}$ form the edges of $G$.

(5) Only those pairs of vertices $v_i, v_j$ are available for contraction for which the principle sub-matrix formed by the elements in the intersection of $i$-th row/column and $j$-th row/column has the form $\begin{bmatrix} 0 & 1 \\ 1 & 0 \end{bmatrix}$.

(6) All those principle sub-matrices formed by elements in the intersection of rows/columns having labels $i_1, i_2, \cdots, i_r$ such that $a_{ii} = 0$, and $a_{ij} = k$ such that $k \geq 2$ for all $i, j \in \{i_1, i_2, \cdots, i_r\}$, then the set of vertices $\{v_{i_1}, v_{i_2}, \cdots, v_{i_r}\}$ forms a partite set or an independent set forming a color class.

(7) Only the diagonal elements of $T(G)$ are zero and all the other elements of $T(G)$ are greater than zero.

(8) In one contraction obtained by identifying some two adjacent vertices, i.e. vertices at unit distance, the size of $T(G)$ reduces by one unit and the count of edges reduces from say $q$ to $q-$ (number of triangles (complete graphs on three points) having the contracted edge as one of their edges + 1).

(9) The collection of subsets of vertices forming the principle sub-matrices of $T(G)$ containing off-diagonal elements strictly bigger than one, and covering all the vertices as a member of some of these subsets, forms partite-representation for $G$. When the count of subsets forming partite-representation is minimal then it is called a minimal-partite-representation.

(10) Every edge in the graph has end points belonging to different partite sets.

(11) When $G$ is $k$-chromatic there will exist $k$ number of principle sub-matrices (that forms the minimal-partite-representation) containing all the off-diagonal elements strictly greater than one, and each such a matrix made up of rows/columns corresponding to certain set of vertices forming a partite set, such that every vertex of $G$ gets incorporated in some of these disjoint partite sets.

(12) The order of largest principle sub-matrix having all off-diagonal elements equal to unit is the order of the largest complete subgraph of $G$.

(13) The order of largest principle sub-matrix having all off-diagonal elements strictly greater than unit represents the independence number of graph $G$.

(14) The nonexistence of at least $k$ principle sub-matrices having all the off-diagonal entries strictly bigger than one implies less than $k$-colorable nature of the graph under consideration.

(15) If removal of an edge or vertex reduces the count of the partite sets, for the minimal-partite-representation, from $k$ to $l$, $l < k$, then the graph under consideration is $k$-critical.

(16) A row (column) of $T(G)$, like the adjacency matrix, $A(G)$, of every connected graph $G$ contains $d$ number of units where $d$ is the degree of the vertex.

(17) If there is not a single principle sub-matrix of $T(G)$ of size $\geq 2$ containing elements $\geq 2$ then the partite sets are singletons (principle sub-matrix of $T(G)$ of size one and containing element zero).

**5. Contraction Algorithm:** We note the effect of a contraction on the transparency matrix and then proceed to develop contraction algorithm which will take a $k$-chromatic graph to a complete graph of largest possible size. When $v_i, v_j$ are adjacent vertices then by $(v_i \Rightarrow v_j)$ we denote the contracting of edge $(v_i, v_j)$ and identifying the vertex $v_i$ with the vertex $v_j$. Let $\tilde{G}$ be the graph that results after the operation $(v_i \Rightarrow v_j)$ on $G$, and let $T(\tilde{G})$ denotes the transparency matrix for $\tilde{G}$, then

$$T(\tilde{G}) = [a_{ij}]_{(p-1)\times(p-1)}$$

can be obtained from $T(G)$ by performing the following operations:
(1) Replace all elements $a_{jk}$, $k \neq j$ by min $\{a_{ik}, a_{jk}\}$.
(2) Replace all elements $a_{kj}$, $k \neq j$ by min $\{a_{ki}, a_{kj}\}$.
(3) Delete $i$-th row and $i$-th column.
(4) If the shortest path (deciding distance) joining vertices $v_m$ and $v_n$ contains the edge $(v_i, v_j)$ then replace the entries $a_{mn}$ and $a_{nm}$ by respectively $(a_{mn} - 1)$ and $(a_{nm} - 1)$.
(5) Keep all other elements as they are.

**Remark 5.1:** How one finds out whether the condition mentioned in (4) is true or not? The answer is simple: The condition will be true if and only if $a_{mn} = a_{mi} + a_{jn} + 1$. Thus, in case of a tree graph the condition will be true when the vertices $v_m$ and $v_n$ belong to different connected components when the branch represented by edge $(v_i, v_j)$ is removed.

**Example 5.1:** Let the graph $G$ be a 5-cycle:
$$v_1 \to v_2 \to \cdots \to v_5 \to v_1,$$
clearly,

$$T(G) = \begin{bmatrix} 0 & 1 & 2 & 2 & 1 \\ 1 & 0 & 1 & 2 & 2 \\ 2 & 1 & 0 & 1 & 2 \\ 2 & 2 & 1 & 0 & 1 \\ 1 & 2 & 2 & 1 & 0 \end{bmatrix}$$

Partite sets forming minimal-partite-representation are (look at the corresponding matrices) {1}, {2, 4}, {3, 5}. Let $\tilde{G}$ be the graph that results after contraction operation ($v_1 \Rightarrow v_2$) which will reduce the graph to a 4-cycle: $v_2 \to v_3 \to v_4 \to v_5 \to v_2$ and will lead to

$$T(\tilde{G}) = \begin{bmatrix} 0 & 1 & 2 & 1 \\ 1 & 0 & 1 & 2 \\ 2 & 1 & 0 & 1 \\ 1 & 2 & 1 & 0 \end{bmatrix}$$

Partite sets forming minimal-partite-representation will be now {2, 4}, {3,5}. If we will carry out further contraction operation ($v_2 \Rightarrow v_3$) we will get the new graph isomorphic to $K_3$, a complete graph on three points, having partite sets: {3}, {4}, {5}.

We now proceed with the contraction algorithm:

**Algorithm 5.1:**
(1) Construct transparency matrix, $T(G)$, for the given $k$-chromatic graph, $G$, on $p$ points.
(2) Find the pair of rows with label $m$ and $n$ such that $a_{mn} = 1$ and one row between the rows contains maximum number of units and the other row contains minimum number of units (or maximum number of non-units), such that **maximum number of replacements** of non-units by units will take place in the contraction operation.
(3) Carry out contraction operation ($v_m \Rightarrow v_n$) resulting in a graph on $(p-1)$ points, say $G_1$.

(4) Construct transparency matrix for the graph, $G_1$, that results after this contraction operation, say $T(G_1)$.

(5) Rename $G_1$ as $G$ and go to step (1) till all the off-diagonal elements of $T(G_1)$ for the graph $G_1$ that results after a contraction operation become units. □

It is clear that the above discussed algorithm will **definitely terminate in finitely many steps**, terminating into formation of a complete graph. The algorithm consists of a sequence of contraction operations in a **preferred way** (causing replacement of maximum number of non-units by units) so that these contractions will lead to formation of complete graph of largest possible size (the size of the complete graph that results will be largest if it is achieved in minimum number of steps, i.e. minimum number of contraction operations).

**Example 5.2:**

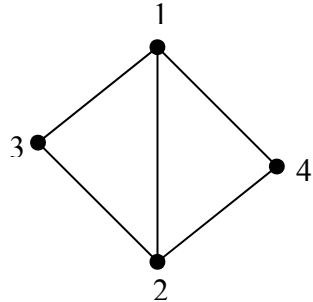

$$G$$

$$T(G) = \begin{bmatrix} 0 & 1 & 1 & 1 \\ 1 & 0 & 1 & 1 \\ 1 & 1 & 0 & 2 \\ 1 & 1 & 2 & 0 \end{bmatrix}$$

Partite sets are: $\{1\}, \{2\}, \{3, 4\}$.
Note that $(1 \Rightarrow 2)$ is **not correct** contraction operation, while, in fact, any other contraction operation is **correct** as well as equivalent as per the above algorithm 5.1. We carry out $(1 \Rightarrow 3)$ which will produce $K_3$, a complete graph on three points having partite sets: $\{2\}, \{3\}, \{4\}$.

Now, given a *k*-chromatic graph on $p$ points and $q$ lines then how many contractions are required to reduce it to a complete graph? In this respect we give the following result.

**Remark 5.2:** As per the **property (8)** among the properties of $T(G)$ given above we know that every contraction reduces the size of the transparency matrix by **exactly one count**. In the light of this property we now proceed with the Hadwiger's conjecture.

Let $G$ be a graph on $p$ points, $p = l + k$. If $G$ reduces to a complete graph (i.e. every non diagonal element of $T(G^*)$ becoming unity where $G^*$ is the graph that results from $G$ after contractions) in $\leq l$ contractions then $G$ will be $\geq k$ chromatic. On the contrary, if $G$ is *k*-chromatic and if all possible $l$ contractions carried out in all possible orders does not produce graph $G^*$ such that the nondiagonal element of its $T(G^*)$ are units then we need to see that the original graph $G$ cannot be *k*-chromatic.

**Theorem 5.1:** Let $G$ be a graph on $p$ points, $p = l + k$, which is either itself *k*-chromatic or a less than *k*-chromatic but the one that has arrived at by suitable contractions in the sense of algorithm 5.1 from initially a *k*-chromatic graph then $G$ is contractible to $K_k$, a complete graph on $k$ points, in at most $k$ contractions.

**Proof:** We proceed by induction on $l$.
**Step 1:** (i) $l = 0$. In this case $G \cong K_k$ and the case is clear.
(ii) $l = 1$. In this case if $G$ is *k*-chromatic then it is isomorphic to graphs which must contain $K_k$ as a subgraph. Or, if $G$ is less than *k*-chromatic but has arrived at by suitable contractions from a *k*-chromatic graph then it must contain a graph isomorphic to graphs like $K_k - x$ or others as a subgraph which go to $K_k$ by only one further suitable contraction in the sense of algorithm 5.1, where $x$ is some edge of $K_k$, otherwise, $G$ cannot be *k*-chromatic.
(iii) $l = 2$. In this case $G$ is isomorphic to graphs each one of which contain one or more suitable edges (in the sense of the algorithm 5.1), one of which is to be contracted, and when it will be contracted this contraction will produce a graph among the graphs suitable for the case

$l = 1$.

**Step 2:** Let the claim be true for all integers $< l$, i.e. the *k*-chromatic graphs or graphs which have arrived by suitable contractions from a *k*-chromatic graph but are themselves less than *k*-chromatic and containing points $p + s$, where $0 \leq s \leq (l-1)$ are only those graphs which contain a contraction, suitable in the sense of algorithm 5.1, which when carried out produces a graph on points $p + (s-1)$, again among the graphs suitable in the sense of algorithm 5.1, so that the process of contraction can be furthered up to $K_k$.

We now proceed to see that it should be valid for $l$. But this is clear since graphs on *p* points, $p = l + k$, are either themselves *k*-chromatic, or are graphs which themselves could be less than *k*-chromatic but have arrived at by suitable contractions, in the sense of algorithm 5.1, from a graph which was *k*-chromatic. Such graphs must be those which must contain a suitable contraction in the sense of algorithm 5.1 which when carried out produce graphs isomorphic to some graph among the suitable graphs (to be *k*-chromatic) for the case with $p = k + (l-1)$, for if these graphs do not contain a contraction which will produce a graph isomorphic to some graph among the graphs which are suitable (for *k*-chromaticity) for the case of $(l-1)$ then $G$ cannot be *k*-chromatic. □

**Theorem 5.2 (Hadwiger's Conjecture):** Every *k*-chromatic graph is contractible to $K_k$, a complete graph on *k* points.

**Proof:** Straightforward from theorem 5.1 given above. □